# VARIATION AND OSCILLATION INEQUALITIES FOR CONVOLUTION PRODUCTS

KARIN REINHOLD AND ANNA K. SAVVOPOULOU

ABSTRACT. We establish variation and oscillation inequalities for convolution products of probability measures on $\mathbb{Z}$.

## 1. Preliminaries

1.1. **Introduction and definitions.** Oscillation and variation inequalities were first used in the ergodic theory context by Bourgain [3] as a means of establishing a.e convergence on a dense set of functions. These results were extended in [4]. The study of oscillation and variation inequalities for operators $\mu^n f$ (see definition below) where $\mu^n$ is the $n$th convolution of a probability measure on a locally compact abelian group were studied in [5]. In this paper we study oscillation and variation inequalities for operators $\mu_n f$ where $\mu_n$ is the convolution of not necessarily distinct probability measures on $\mathbb{Z}$.

**Definition 1.1.** Given a set of real numbers $\{x_n\}_{n \in I}$, where $I$ is a countable index set, define its variation $\rho$ norm by

$$\|x_n\|_{v(\rho)} = \sup_{n(j)} \left( \sum_{j=1}^{\infty} |x_{n(j)} - x_{n(j+1)}|^{\rho} \right)^{\frac{1}{\rho}}$$

where the supremum is taken over all possible sequences $n(1) \leq n(2) \leq \cdots$ and each $n(j) \in I$.

Unless otherwise noted, we assume that $I = \mathbb{Z}^+$. In case $I = \{n : n_k \leq n < n_{k+1}\}$ for a subsequence $n_k$, we will write $\|x_n : n_k \leq n < n_{k+1}\|_{v(\rho)}$. We also define the oscillation norm of a sequence $\{x_n\}_{n \in I}$.

**Definition 1.2.** Given an increasing sequence of positive integers $\{n_k\}_{k=1}^{\infty}$, define the oscillation norm of a sequence of numbers $\{x_n\}_{n=1}^{\infty}$ by

$$\|x_n\|_{o(s)} \left( \sum_{k=1}^{\infty} \max_{n_k \leq n \leq n_{k+1}} |x_n - x_{n_k}|^s \right)^{1/s}.$$

We will be studying the variation and oscillation of operators given by sequences of convolution products of measures on $\mathbb{Z}$, where the convolution of two probability measures $\nu_1$ and $\nu_2$ on $\mathbb{Z}$ is given by

$$\nu_1 * \nu_2(k) = \sum_{j \in \mathbb{Z}} \nu_1(j)\nu_2(k - j).$$

**Definition 1.3.** Let $(X, \Sigma, m)$ be a probability measure space and $T : X \to X$ be an invertible measure preserving transformation, i.e $m(T^{-1}B) = m(B) \, \forall B \in \Sigma$. If $\{\nu_i\}_{i=1}^{\infty}$ is a sequence of probability measures on $\mathbb{Z}$ we let $\mu_n = \nu_1 * \cdots * \nu_n$ and we define $\mu_n f$ for $f \in L^1(X)$ by

$$\mu_n f(x) = \sum_{j \in \mathbb{Z}} \mu_n(j) f(T^j x).$$

A few more definitions are essential before we state our main result.

**Definition 1.4.** A probability measure $\mu$ defined on a group $G$ is called strictly aperiodic if and only if the support of $\mu$ can not be contained in a proper left coset of $G$.

**Definition 1.5.** If $p > 0$, the $p^{\text{th}}$ moment of $\mu$ is given by $\sum_{k \in \mathbb{Z}} |k|^p \mu(k)$ and it is denoted by $m_p(\mu)$. The expectation of $\mu$ is $\sum_{k \in \mathbb{Z}} k \mu(k)$ and is denoted by $E(\mu)$.

**Notation 1.6.** Throughout the paper $\|\mu_n f\|_{v(\rho)}$ or $\|\mu_n f\|_{o(\rho)}$ will be understood to stand for $\|\mu_n f(x)\|_{v(\rho)}$ or $\|\mu_n f(x)\|_{o(\rho)}$ where $x$ is some fixed member of $X$.

Our main result is the following.

**Theorem 1.7.** Let $(\nu_n)$ be a sequence of strictly aperiodic probability measures on $\mathbb{Z}$ such that
  (1) $E(\nu_n) = 0$, $\forall n$
  (2) $\sum_{i=1}^{n} m_2^2(\nu_i) = O(n)$
  (3) there exist a constant $C$ and an integer $N_0 > 0$, such that $|\hat{\nu}_n(t)| \leq e^{-Ct^2}$, $\forall n > N_0$ and $t \in [-1/2, 1/2)$.

Let $\mu_n = \nu_1 * \cdots * \nu_n$, then for $2 < s < \infty$ the following inequalities hold
(1) $$\|\|\mu_n f\|_{v(s)}\|_2 \leq C \|f\|_2 \text{ and}$$
(2) $$\|\|\mu_n f\|_{o(2)}\|_2 \leq C \|f\|_2$$

**Remark 1.8.** If we assume that the second moments are uniformly bounded then condition 3 follows from Theorem 3.3 of [6].

**Remark 1.9.** Let

$$\nu_n(k) = \begin{cases} \dfrac{1 - a_n}{2} & k = \pm 1 \\ a_n & k = 0 \\ 0 & \text{otherwise} \end{cases}$$

where $1 > a_n > 0$ and $a_n \to 0$ fast enough so that $\prod_{n=1}^{\infty} a_n > 0$. Then, using an argument similar to that in [1] one may show that the sequence $\mu_n f$ does not converge a.e for some $f \in L^\infty$. Of course the sequence $\nu_n(k)$ does not satisfy the condition $\sup_n \sup_{\alpha,\beta} \nu_n(\beta \mathbb{Z} + \alpha) \leq \rho$ while it does satisfy the condition $m_1(\nu_n) \leq a$.

## 1.2. Properties of variation norms.

**Proposition 1.2.1** (Properties of the $v(\rho)$−norms).
(1) *For each $\rho$, $1 \leq \rho < \infty$, $\|\cdot\|_{v(\rho)}$ is a semi-norm.*
(2) *If $\{x_n\}$ is a sequence such that $x_{n_k} = 0$ for each $k$, then*
$$\|x_n\|_{v(\rho)} \leq 2 \left( \sum_k \|x_n : n_k \leq n < n_{k+1}\|_{v(\rho)}^\rho \right)^{1/\rho}.$$
(3) $\|x_n\|_{v(\rho)} \leq 2 \left( \sum_n |x_n|^\rho \right)^{1/\rho}.$

For a proof see [4].

## 1.3. Proposition 1.3.1 and proof of (3) in Proposition 1.3.1.

As in [5] we establish our main result by considering two related inequalities. Proposition 1.3.1 ties these inequalities together to produce the result of Theorem 1.7. This section deals with the proof of condition (3) of Proposition 1.3.1. Condition (3) is established in Theorem 1.3.3 and it relates $\mu_{4^k} f$ to its symmetrized version $\lambda_{4^k} f$, where $\lambda_n = \nu_1 * \cdots * \nu_{4^k} * \tilde{\nu}_1 * \cdots * \tilde{\nu}_{4^k}$ and $\tilde{\nu}_i(k) = \nu_i(-k)$.

**Proposition 1.3.1.** *Suppose that the following inequalities hold:*

$$(3) \qquad \left\| \left( \sum_{k=1}^\infty |\mu_{4^k} f - \lambda_{4^k} f|^2 \right)^{1/2} \right\|_2 \leq C \|f\|_2$$

$$(4) \qquad \left\| \sum_{k=1}^\infty \|\mu_n f : 4^{k-1} \leq n < 4^k\|_{v(2)}^2 \right\|_2 \leq C \|f\|_2$$

*then*
$$\|\|\mu_n f\|_{v(\rho)}\|_2 \leq C \|f\|_2 \; \forall \rho > 2 \text{ and}$$
$$\|\|\mu_n f\|_{o(2)}\|_2 \leq C \|f\|_2.$$

*Proof.* The proof is similar to the proof of Theorem 4.4 in [5]. □

The inequalities in the assumption of Proposition 1.3.1 are proved in Theorems 1.3.3 and 1.4.1. The first inequality relates $\mu_n f$ to a symmetrized convolution product $\lambda_n f$, while the second inequality examines the variation of $\mu_n f$ over finite intervals of integers.

**Lemma 1.3.2.** *Let $(\nu_n)$ be a sequence of strictly aperiodic probability measures on $\mathbb{Z}$ with finite second moment such that*

(1) $E(\nu_n) = 0$ , $\forall n$
(2) *There exist a constant $C$ and an integer $N_0 > 0$, such that $|\hat{\nu}_n(t)| \leq e^{-Ct^2}$, $\forall n > N_0$ and $t \in [-1/2, 1/2)$ .*

*Let $\mu_n = \nu_1 * \cdots * \nu_n$, then*

$$|1 - \hat{\mu}_n(t)| \leq t^2 \sum_{i=1}^{n} a_i$$

*where $a_i = cm_2(\nu_i)$ for some uniform constant $c$.*

*Proof.*

$$\begin{aligned}
|1 - \hat{\mu}_n(t)| &= \left| -\int_0^t (\hat{\mu}_n(s))' \, ds \right| \\
&= \left| -\int_0^t \left( \prod_{i=1}^{n} \hat{\nu}_i(s) \right)' \, ds \right| \\
&= \left| \int_0^t \sum_{i=1}^{n} \prod_{\substack{j=1 \\ j \neq i}}^{n} \hat{\nu}_j(s) \hat{\nu}_i'(s) \, ds \right| \\
&\leq \sum_{i=1}^{n} \int_0^t \prod_{\substack{j=1 \\ j \neq i}}^{n} |\hat{\nu}_j(s)| |\hat{\nu}_i'(s)| \, ds \\
&\leq \sum_{i=1}^{n} \int_0^t a_i |s| \, ds \\
&\leq \sum_{i=1}^{n} a_i t^2
\end{aligned}$$

□

**Theorem 1.3.3.** *Let $(\nu_n)$ be a sequence of strictly aperiodic probability measures on $\mathbb{Z}$ such that*

(1) $E(\nu_n) = 0$ , $\forall n$
(2) $\phi(n) = \sum_{i=1}^{n} m_2(\nu_i) = O(n)$

(3) There exist a constant $C$ and an integer $N_0 > 0$, such that $|\hat{\nu}_n(t)| \leq e^{-Ct^2}$, $\forall n > N_0$ and $t \in [-1/2, 1/2)$

and let $\lambda = \mu * \tilde{\mu}$, where $\tilde{\mu}(x) = \mu(-x)$ for $x \in \mathbb{Z}$, then

$$\left\| \left( \sum_{k=1}^{\infty} |\mu_{4^k} f - \lambda_{4^k} f|^2 \right)^{\frac{1}{2}} \right\|_2 \leq C \|f\|_2$$

*Proof.* Note that

$$\hat{\tilde{\mu}}_n(t) = \hat{\mu}_n(-t) = \sum_{k \in \mathbb{Z}} \mu_n(k) e^{2\pi i k t} = \bar{\hat{\mu}}_n(t), \text{ for } t \in [-1/2, 1/2)$$

thus

$$|1 - \hat{\tilde{\mu}}_n(t)| = |1 - \hat{\mu}_n(t)|.$$

A standard technique known as the Calderón Transfer Principle ( [2]) allows us to establish

$$\left\| \left( \sum_{k=1}^{\infty} |(\mu_{4^k} - \lambda_{4^k}) f|^2 \right)^{\frac{1}{2}} \right\|_2^2 \leq K \|f\|_2^2$$

for all $f \in L^2(X)$ and for some constant $K$, by showing that

$$\left\| \left( \sum_{k=1}^{\infty} |(\mu_{4^k} - \lambda_{4^k}) * f|^2 \right)^{\frac{1}{2}} \right\|_2^2 \leq K \|f\|_2^2$$

for all $f \in l_2(\mathbb{Z})$ for the same constant $K$.

By Parsevals' identity

$$\left\| \left( \sum_{k=1}^{\infty} |\mu_{4^k} * f - \lambda_{4^k} * f|^2 \right)^{\frac{1}{2}} \right\|_2^2 = \int_{-1/2}^{1/2} \sum_{k=1}^{\infty} \left| \hat{\mu}_{4^k}(t) - \hat{\lambda}_{4^k}(t) \right|^2 |\hat{f}(t)|^2 \, dt$$

Note that by the assumptions and Lemma 1.3.2 for $\rho > 1$ we have that

$$
\begin{aligned}
\sum_{k=1}^{\infty} |\hat{\mu}_{\rho^k}(t) - \hat{\lambda}_{\rho^k}(t)|^2 &= \sum_{k=1}^{\infty} \left| \hat{\mu}_{\rho^k}(t) - \hat{\mu}_{\rho^k}(t)\hat{\tilde{\mu}}_{\rho^k}(t) \right|^2 \\
&\leq t^4 \sum_{k=1}^{\infty} e^{-2c\rho^k t^2} \left( \sum_{i=1}^{\rho^k} a_i \right)^2 \\
&\leq ct^4 \sum_{k=1}^{\infty} \rho^{2k} e^{-2c\rho^k t^2} \\
&\leq ct^4 \int_1^{\infty} \rho^{2x} e^{-2c\rho^x t^2} \, dx \\
&= \frac{Ct^4}{t^2 \ln(\rho)} \int_{2c\rho t^2}^{\infty} \frac{u}{2t^2} e^{-u} \, du \\
&< \frac{C}{\ln(\rho)},
\end{aligned}
$$

and the claim follows for some constant $K$. □

### 1.4. Proof of (4) in Proposition 1.3.1. Conclusion of proof, of main theorem.

**Theorem 1.4.1.** *Let $(\nu_n)$ be a sequence of strictly aperiodic probability measures on $\mathbb{Z}$ such that*

(1) $E(\nu_n) = 0$, $\forall n$
(2) $\sum_{i=1}^{n} m_2^2(\nu_i) = O(n)$
(3) *There exist a constant $C$ and an integer $N_0 > 0$, such that $|\hat{\nu}_n(t)| \leq e^{-Ct^2}$, $\forall n > N_0$ and $t \in [-1/2, 1/2)$*

*and $\mu_n = \nu_1 * \cdots * \nu_n$, then*

$$\left\| \left( \sum_{k=1}^{\infty} \|\mu_n f : 4^{k-1} \leq n < 4^k\|_{v(2)}^2 \right)^{1/2} \right\|_2 \leq C\|f\|_2$$

Following notation in [4] we will denote by $\mathcal{S}$ the two norm of a sequence, i.e

$$\mathcal{S}(x_n) = \left( \sum_{n=1}^{\infty} |x_n|^2 \right)^{1/2}.$$

Additionally we will let $K_n = \mu_n - \mu_{4^{k-1}}$ for $4^{k-1} \leq n < 4^k$. With this notation we are trying to prove

(5) $$\|\mathcal{S}(\|K_n f : 4^{k-1} \leq n < 4^k\|_{v(2)})\|_2 \leq C\|f\|_2\,.$$

**Definition 1.4.2.** For $j \geq 0$ let $E_j = [-2^{-j}, -2^{-j-1}) \cup [2^{-j-1}, 2^{-j})$ and $\chi_j(t) = \chi_{E_j}(\{\text{sgn}(t)\}(t - \lfloor t \rfloor))$, where $\chi_j(t)$ is seen to be the one-periodic extension of $\chi_{E_j}$. For $j \leq 0$ let $\chi_j(t) = 0$.

**Definition 1.4.3.** With $E_j$ as in Definition 1.4.2 and $4^{k-1} \leq n < 4^k$ let $K_{j,n}$ denote the sequence whose Fourier transform satisfies

$$\hat{K}_{j,n}(t) = \hat{K}_n(t)\chi_{j+k}(t)\,.$$

**Remark 1.4.4.** By the definition of $K_{j,n}$ we have that $K_n = \sum_{j \in \mathbb{Z}} K_{j,n}$.

**Definition 1.4.5.** Define the norm $\|\cdot\|$ by $\|\cdot\| = \|\|\cdot\|_{v(2)}\|_2$.

**Lemma 1.4.6.** *The operator $\mathcal{S}$ commutes with the $L^2$-norm and therefore inequality (2) can be written as $\mathcal{S}\left(\|K_n f : 4^{k-1} \leq n < 4^k\|\right) \leq C\|f\|_2$. Furthermore the above inequality is true if*

$$\mathcal{S}\left(\|K_{j,n} f : 4^{k-1} \leq n < 4^k\|\right) \leq \frac{C}{4^{|j|/2}}\|f\|_2$$

*Proof.* To show that $\mathcal{S}$ commutes with the two norm, we observe that

$$\begin{aligned}
\|\mathcal{S}(\mu_n f)\|_2 &= \left(\int \sum_{n=1}^{\infty} |\mu_n f(x)|^2 \, dx\right)^{1/2} \\
&= \left(\sum_{n=1}^{\infty} \int |\mu_n f(x)|^2 \, dx\right)^{1/2} \\
&= \left(\sum_{n=1}^{\infty} \|\mu_n f\|_2^2\right)^{1/2} \\
&= \mathcal{S}(\|\mu_n f\|_2).
\end{aligned}$$

Therefore setting $x_k = \|K_n f : 4^{k-1} \leq n < 4^k\|_{v(2)}$ we have that

$$\|\mathcal{S}(x_k)\|_2 = \mathcal{S}(\|x_k\|_2) = \mathcal{S}(\|K_n f : 4^{k-1} \leq n < 4^k\|)\,.$$

Now suppose $\mathcal{S}\left(\|K_{j,n}f : 4^{k-1} \leq n < 4^k\|\right) \leq \frac{C}{4^{|j|/2}}\|f\|_2$. Then

$$\begin{aligned}
\mathcal{S}\left(\|K_n f : 4^{k-1} \leq n < 4^k\|\right) &\leq \sum_j \mathcal{S}\left(\|K_{j,n}f : 4^{k-1} \leq n < 4^k\|\right) \\
&\leq \sum_j \frac{C}{4^{|j|/2}}\|f\|_2 \\
&\leq C\|f\|_2.
\end{aligned}$$

$\square$

Thus it suffices to prove that

$$\mathcal{S}\left(\|K_{j,n} : 4^{k-1} \leq n < 4^k\|\right) \leq \frac{C}{4^{|j|/2}}\|f\|_2.$$

**Lemma 1.4.7.**
(1) *For every $j$ and $n$ we have*

$$|\hat{K}_{j,n}(t)| \leq \frac{C}{4^{|j|}}$$

(2) *For every $j$ and $n$, if $4^{k-1} \leq n < 4^k$, we have*

$$|\hat{K}_{j,n}(t) - \hat{K}_{j,n+1}(t)| \leq \frac{C}{4^k} m_2(\nu_{n+1})$$

*Proof.*

(1) Note that $\hat{K}_{j,n}(t) \neq 0$ if and only if $2^{-j-k-1} \leq |t| \leq 2^{-j-k}$. We will first show that
$|\hat{K}_n(t)| \leq C 4^k t^2$. Observe that for $4^{k-1} \leq n < 4^k$

$$\begin{aligned}
|\hat{K}_n(t)| &= |\hat{\mu}_n(t) - \hat{\mu}_{4^{k-1}}(t)| \\
&= |\hat{\mu}_{4^{k-1}}(t)||1 - \prod_{i=4^{k-1}+1}^{4^k} \hat{\nu}_i(t)| \\
&\leq t^2 \sum_{i=4^{k-1}+1}^{n} m_2(\nu_i) \\
&\leq t^2 B 4^k
\end{aligned}$$

then $|\hat{K}_{j,n}(t)| \leq B4^k 2^{-2j-2k} = B\dfrac{1}{4^j} = \dfrac{B}{4^{|j|}}$ for positive $j$.

We also have
$$\begin{aligned}|\hat{K}_n(t)| &\leq 2e^{-4^{k-1}Ct^2} \\ &\leq \dfrac{2}{(C/4)4^k t^2} \\ &= \dfrac{8}{C4^k t^2}.\end{aligned}$$

Then
$$|\hat{K}_{j,n}(t)| \leq \dfrac{4}{C4^k t^2} \leq 8\dfrac{2^{2(j+k+1)}}{C4^k} = C'4^j \quad \text{for negative } j$$
so that
$$|\hat{K}_{j,n}(t)| \leq C\dfrac{1}{4^{|j|}}.$$

(2) If $t \in E_{j+k}$ and $4^{k-1} \leq n < n+1 \leq 4^k$
$$\begin{aligned}|\hat{K}_{j,n}(t) - \hat{K}_{j,n+1}(t)| &= |\hat{\mu}_n(t) - \hat{\mu}_{n+1}(t)| = |\hat{\mu}_n(t)||\hat{\nu}_{n+1}(t) - 1| \\ &\leq Ce^{-nCt^2} a_{n+1} t^2 \\ &\leq C'\dfrac{m_2(\nu_{n+1})}{4^k}. \quad \square\end{aligned}$$

Let $k \in \mathbb{Z}^+$ and $j \in \mathbb{Z}$. We will let
$$I_{j,k} = \{t : 4^{-j-k-1} \leq |t| \leq 4^{-j-k}\},$$
and the points
$$4^{k-1} = \alpha_1 < \alpha_2 < \cdots < \alpha_N < 4^k$$
be $N = 3\min\{4^{k-1}, 4^{|j|}\}$ equidistributed points.

**Definition 1.4.8.** With the notation above define the operators $A_n$ by $A_n = K_{j,\alpha_{m-1}}$ where $\alpha_{m-1} \leq n < \alpha_m$.

**Lemma 1.4.9.** $\mathcal{S}(\|A_n f : 4^{k-1} \leq n < 4^k\|) \leq \dfrac{C}{4^{|j|/2}}\|f\|_2$

*Proof.* For fixed $k$ using the third property as described in Proposition 1.2.1 for $\rho = 2$ we have that
$$\begin{aligned}\|A_n f : 4^{k-1} \leq n < 4^k\| &= \|K_{j,\alpha_m} f : m \leq N\| \\ &\leq 2\left(\sum_{m \leq N} \|K_{j,\alpha_m} f\|_2^2\right)^{1/2}\end{aligned}$$

By Parseval's formula, we have that

$$\sum_{m \leq N} \|K_{j,\alpha_m} f\|_2^2 = \sum_{m \leq N} \int_{[-1/2,1/2)} |\hat{K}_{j,\alpha_m}(t)\hat{f}(t)|^2 \, dt,$$

using Lemma 1.4.7 and noting that $\hat{K}_{j,\alpha_m}(t)$ is non-zero only if $t \in I_{j,k}$, the above is

$$\leq C \sum_{m \leq N} \frac{1}{4^{2|j|}} \int_{I_{j,k}} |\hat{f}(t)|^2 \, dt$$

$$\leq \frac{C'}{4^{|j|}} \int_{I_{j,k}} |\hat{f}(t)|^2 \, dt.$$

In the last step we have used the fact that $N \leq 3 \cdot 4^{|j|}$. Summing the above estimate in $k$ we obtain

$$\sum_k \|A_n f : 4^{k-1} \leq n < 4^k\|^2 \leq \frac{C}{4^{|j|}} \sum_k \int_{I_{j,k}} |\hat{f}(t)|^2 \, dt$$

$$\leq \frac{C}{4^{|j|}} \int_{[-1/2,1/2)} |\hat{f}(t)|^2 \, dt,$$

which by Parseval's formula, implies the following bound

$$\mathcal{S}(\|A_n f : 4^{k-1} \leq n < 4^k\|) \leq \frac{C}{4^{|j|/2}} \|f\|_2.$$

$\square$

**Lemma 1.4.10.** $\mathcal{S}(\|K_{j,n} f : 4^{k-1} \leq n < 4^k\|) \leq \frac{C}{4^{|j|/2}} \|f\|_2$

*Proof.* We have

$$\mathcal{S}(\|K_{j,n} f : 4^{k-1} \leq n < 4^k\|) \leq$$
$$\mathcal{S}(\|A_n f : 4^{k-1} \leq n < 4^k\|) + \mathcal{S}(\|A_n f - K_{j,n} f : 4^{k-1} \leq n < 4^k\|).$$

We need only prove that $\mathcal{S}(\|A_n f - K_{j,n} f : 4^{k-1} \leq n < 4^k\|) \leq \frac{C}{4^{|j|/2}} \|f\|_2$.

Since $A_n = K_{j,n}$, if $4^{k-1} \leq n < 4^k$ and $k \leq |j|$, we assume that $k > |j|$ and hence $N = 3 \cdot 4^{|j|}$.

$$\|A_n f - K_{j,n} f : 4^{k-1} \leq n < 4^k\|_{v(2)} \leq$$
$$2\mathcal{S}(\|A_n f - K_{j,n} f : \alpha_{m-1} \leq n < \alpha_m\|_{v(2)}) =$$
$$2\mathcal{S}(\|K_{j,n} f : \alpha_{m-1} \leq n < \alpha_m\|_{v(2)}) =$$
$$2 \left( \sum_{m \leq N} \|K_{j,n} f : \alpha_{m-1} \leq n < \alpha_m\|_{v(2)}^2 \right)^{1/2}.$$

For fixed $m$ we have

$$\|K_{j,n}f : \alpha_{m-1} \leq n < \alpha_m\|_{v(2)} \leq$$
$$\|K_{j,n}f : \alpha_{m-1} \leq n < \alpha_m\|_{v(1)} \leq$$
$$\sum_{\alpha_{m-1} \leq n < \alpha_m} |K_{j,n}f - K_{j,n+1}f| \leq$$
$$(\alpha_m - \alpha_{m-1})^{1/2} \left( \sum_{\alpha_{m-1} \leq n < \alpha_m} (K_{j,n}f - K_{j,n+1}f)^2 \right)^{1/2}.$$

If we put the above two estimates together, take the 2−norm and use that $S$ commutes with the 2−norm, we get

$$\|A_n f - K_{j,n}f : 4^{k-1} \leq n < 4^k\| =$$
$$\| \|A_n f - K_{j,n} : 4^{k-1} \leq n < 4^k\|_{v(2)} \|_2 \leq$$
$$2S \left( (\alpha_m - \alpha_{m-1})^{1/2} \left( \sum_{\alpha_{m-1} \leq n < \alpha_m} \|K_{j,n}f - K_{j,n+1}f\|_2^2 \right)^{1/2} \right) =$$
$$2 \left( \sum_{m \leq N} (\alpha_m - \alpha_{m-1}) \left( \sum_{\alpha_{m-1} \leq n < \alpha_m} \|K_{j,n}f - K_{j,n+1}f\|_2^2 \right) \right)^{1/2}.$$

By Parseval's formula and recalling that $\alpha_m - \alpha_{m-1} = 3 \cdot 4^{k-1}/N$, we have

$$\sum_{m \leq N} (\alpha_m - \alpha_{m-1}) \sum_{\alpha_{m-1} \leq n < \alpha_m} \|K_{j,n}f - K_{j,n+1}f\|_2^2$$
$$= \sum_{m \leq N} \frac{3 \cdot 4^{k-1}}{N} \sum_{\alpha_{m-1} \leq n < \alpha_m} \int_{[-1/2,1/2]} |\hat{K}_{j,n}(t) - \hat{K}_{j,n+1}(t)|^2 |\hat{f}(t)|^2 \, dt,$$

by Lemma 1.4.7 and noting that $\hat{K}_{j,n}(t)$, $\hat{K}_{j,n+1}(t)$ is non-zero only if $t \in I_{j,k}$,

$$\leq C \frac{3 \cdot 4^{k-1}}{N} \frac{1}{4^{2k}} \sum_{n \leq 4^k} m_2^2(\nu_{n+1}) \int_{I_{j,k}} |\hat{f}(t)|^2 \, dt$$
$$= \frac{C}{N} \left( \frac{1}{4^k} \sum_{n \leq 4^k} m_2^2(\nu_{n+1}) \right) \int_{I_{j,k}} |\hat{f}(t)|^2 \, dt$$
$$\leq \frac{C'}{N} \int_{I_{j,k}} |\hat{f}(t)|^2 \, dt.$$

Summing the above estimate in $k$ we obtain,

$$\sum_k \|A_n f - K_{j,n} f : 4^{k-1} \leq n < 4^k\|^2 \leq \frac{C}{N} \sum_{k > |j|} \int_{I_{j,k}} |\hat{f}(t)|^2 \, dt$$
$$\leq \frac{C}{3 \cdot 4^{|j|}} \int_{[-1/2, 1/2)} |\hat{f}(t)|^2 \, dt,$$

which, by Parseval's formula, implies the bound

$$\mathcal{S}(\|A_n f - K_{j,n} f : 4^{k-1} \leq n < 4^k\|) \leq \frac{C}{4^{|j|/2}} \|f\|_2$$

and the claim has been proved. $\square$

In summary, inequalities (3) and (4) of Proposition 1.3.1 combine to give Theorem 1.7, the main result of this paper. Inequality (3) is proved by Lemma 1.3.2, while (4) is proved by combining Lemmas 1.4.6, 1.4.7, 1.4.9 and 1.4.10.

(Reinhold) Department of Mathematics, University at Albany, SUNY, Albany, NY 12222
*E-mail address*: `reinhold@albany.edu`

(Savvopoulou) Department of Mathematical Sciences, Indiana University South Bend, South Bend, IN, 46545
*E-mail address*: `annsavvo@iusb.edu`